# Optimal chemotherapy and immunotherapy schedules for a cancer-obesity model with Caputo time fractional derivative

Tuğba Akman Yıldız*[†1], Sadia Arshad [‡2] and Dumitru Baleanu[§3,4]

[1] *Department of Management, University of Turkish Aeronautical Association, 06790 Ankara, Turkey*
[2] *COMSATS Institute of Information Technology, Lahore, Pakistan*
[3] *Department of Mathematics, Çankaya University, 06530 Ankara, Turkey*
[4] *Institute of Space Sciences, Magurele-Bucharest 077125, Romania*

April 16, 2018

## Abstract

This work presents a new mathematical model to depict the effect of obesity on cancerous tumor growth when chemotherapy as well as immunotherapy have been administered. We consider an optimal control problem to destroy the tumor population and minimize the drug dose over a finite time interval. The constraint is a model including tumor cells, immune cells, fat cells, chemotherapeutic and immunotherapeutic drug concentrations with Caputo time fractional derivative. We investigate the existence and stability of the equilibrium points namely, tumor free equilibrium and coexisting equilibrium, analytically. We discretize the cancer-obesity model using L1-method. Simulation results of the proposed model are presented to compare three different treatment strategies: chemotherapy, immunotherapy and their combination. In addition, we investigate the effect of the differentiation order $\alpha$ and the value of the decay rate of the amount of chemotherapeutic drug to the value of the cost functional. We find out the optimal treatment schedule in case of chemotherapy and immunotherapy applied.

---
*Corresponding author
[†]tr.tugba.akman@gmail.com, takman@thk.edu.tr
[‡]sadia_735@yahoo.com
[§]dumitru@cankaya.edu.tr

# 1 Introduction

Obesity is associated with increased incidence and mortality of multiple cancers, with risk ratios correlating directly with body mass index (BMI) in a dose-dependent fashion. Obesity not only affects cancer incidence and long-term cancer-specific mortality, but also impacts on survival and recurrence among those diagnosed with cancer. Multiple variables regulate the obesity-cancer risk equation. Gender and ethnic differences exist; risk ratios for obesity are significantly higher in men than women for colon cancer incidence. For example, Asians suffer increased breast cancer risk at lower BMI compared to non-Asians. Obesity is associated with a state of chronic systemic inflammation. The link between nutrient excess and inflammation is rooted in the chemical nature of nutrients, bioenergetic molecules capable of participating in energy-intensive reactions that are potentially damaging to cells. In addition, obesity is associated with an increased risk of a number of nutrient deficiencies that have been implicated in cancer pathogenesis, including vitamin D, selenium, and magnesium. To understand the influence of obesity on cancer better, Ku-Carrillo et al. developed a mathematical model in [14] that incorporates the interaction between the immune system and adipose cells.

Among evolving treatment modalities, immunotherapy has progressed rapidly into an auspicious field of medicine. Aim of cancer immunotherapy is to assist the immune system to efficiently recognize and attack cancerous cells. This treatment strategy has been investigated in terms of mathematical models, too. For example, an extended mathematical model to depict interactions between cancer cells and adaptive immune system in mouse is proposed in [29] to optimize protocol of immunotherapy with dendritic cell vaccine. Their model includes tumor cells, natural killers, naïve and mature cytotoxic $T$ cells, regulatory $T$ cells, naïve and mature helper $T$ cells, dendritic cells and interleukin 2 cytokine. In the study [5], importance of scheduling immunotherapy is underlined to eradicate and control tumor population. Moreover, immune competition is analyzed in case of tumor population in [6]. In [31], the effect of chemotherapy and immunotherapy on control of the tumor growth has been presented. The combined effect of immunotherapy and chemotherapy is investigated based on optimal control theory in [25].

The aim behind chemotherapy is to eliminate the tumor cells. The cytotoxic drugs, which have a fast dividing rate, are distributed to different parts of the body through the blood circulation system. Chemotherapy also affects mitosis rate of some other kinds of cells that have rapid growth naturally, including hair cells, in an undesirable way [26]. In several clinical contexts, it is important to minimize, or rather, optimize the amount of drug(s) used in order to regulate the potentially lethal side effects of chemotherapy in cancer treatment. Much remarkable work has been done on modeling chemotherapy treatment by using optimal control theory. For example, optimal control is used to examine the ability of a heterogeneous tumor to counteract the chemotherapeutic drug in [22]. To achieve this goal, it is underlined that the drug must be injected at the maximum rate [16]. Models for chemotherapy schedule have been reviewed based on optimal control in [34]. The study [11] offers a discussion of the prospective application of optimal control theory. Moreover, four optimal control problems for chemotherapy schedule are investigated in [11] and different choices of the objective functions in framework of chemotherapy are compared. Combination of chemotherapy and radiotherapy to eradicate the cancer with metastasis is discussed in [12]. An adaptive robust control is proposed to adjust the drug dosages with an extended Kalman filter observer in [32]. An optimal control strategy based on a linear time varying approximation technique is proposed in [13]. On the other hand, an optimal control problem (OCP) with a free final-time is solved for tumor-immune interactions with the aim of minimizing not only the tumor population but also the treatment period in [4]. A model-free method for chemotherapy based on reinforcement learning is proposed in [24] using the closed-loop control. Precisely, they develop an optimal controller using Q-learning algorithm for cancer chemotherapy treatment. In [7], the authors model the chemotherapy treatment based on optimal control theory and the objective is to minimize the tumor cells while keeping the healthy cells above a fixed level. In the study [35], the interactions between tumor cells, immune cells, chemotherapeutic and immunotherapeutic dug concentrations are given by a system of ordinary differential equations. On the other hand, Pillis and Radunskaya have proposed a system for tumor cells, immune cells, normal cells and chemotherapeutic drug concentration in [8]. Recently, effect of obesity on cancer growth has been investigated in [15] by extending the model of Pillis and Radunskaya with fat cells.

Fractional differentiation and integration operators, which are the generalization of classical integer-order counterparts, are capable of capturing memory effects due to their nonlocal nature [28, 33, 18]. It is a useful tool to develop suitable models for describing real-world problems which cannot be expressed by using integer-order differential equations. In the context of tumor-growth models, we can mention the following papers: A cancer model with two immune effectors is investigated in terms of Caputo fractional derivative in [3], while a model for HIV infection of CD4+T cells is generalized with Caputo derivative [30]. Existence and uniqueness of the solution of a cancer model with Caputo-Fabrizio derivative has been justified in [9]. For fractional optimal control problems (FOCPs), we can mention the study [27] where a fractional malaria transmission model is investigated based on optimal control techniques. In addition, HIV infection is investigated in [10] and West Nile virus model with Caputo derivative is presented in [37] in terms of Caputo derivative. On the other hand, an optimal control strategy is proposed for a nonlinear multi-strain tuberculosis model with Caputo derivatives in [36].

The aim of this paper is to investigate a FOCP governed by a Caputo time fractional cancer-obesity model. We note that the proposed model is based on the Sharma and Samanta's model [35], which presents the growth/decay of immune and tumor cells when chemotherapeutic as well as immunotherapeutic drugs have been injected, and Ku-Carrillo et al.'s work [15], which offers a model for cancer-obesity relation. However, influence of obesity to cancer growth is not considered in the model in [35]. Our model expresses the intercommunication of tumor cells, immune cells, fat cells, chemotherapeutic and immunotherapeutic drug concentrations in terms of Caputo time fractional derivative. To the best of our knowledge, this is the first study investigating the optimal treatment strategy for a cancerous tumor growth model in terms of a fractional derivative. Moreover, we investigate the existence and stability of the equilibrium points namely, tumor free equilibrium and coexisting equilibrium, analytically. We solve the optimal control problem using forward-backward sweep method after discretizing the fractional differential equation (FDE) applying the so-called L1 method. We illustrate the contribution of the use of fractional derivatives by solving the FOCP for different orders of differentiation. In addition, we examine the influence of the decay rate of amount of chemotherapeutic drug to the value of the cost functional. To sum up, we find out the optimal treatment schedule in case of chemotherapy and immunotherapy for a generalized cancer-obesity model.

The rest of the paper is organized as follows: In Sec. 2, we introduce the FOCP governed by cancer-obesity model and present the optimality system. In Sec. 3, existence of the equilibrium points and their stability conditions are discussed. In Sec. 4, we explain the discretization of the FOCP. In Sec. 5, we present some numerical results to compare different treatment strategies, namely chemotherapy, immunotherapy and their combination. Then, the paper ends with summary and conclusion.

## 2 Fractional optimal control problem

In the literature, several fractional derivatives have been defined. One of the mostly used fractional differentiation operators is Caputo derivative [28]. In this section, we briefly mention the required definitions and properties of Caputo derivative in order to derive the optimality system.

Firstly, we define the (left) Caputo fractional differentiation operator for $0 < \alpha < 1$ as in the study [28]:

$$^C_A\mathcal{D}^\alpha_t \varphi(t) = \frac{1}{\Gamma(1-\alpha)} \int_A^t \frac{\varphi'(s)}{(t-s)^\alpha}\, ds, \qquad (2.1)$$

The corresponding right differentiation operator is defined as

$$^C_t\mathcal{D}^\alpha_B \varphi(t) = -\frac{1}{\Gamma(1-\alpha)} \int_t^B \frac{\varphi'(s)}{(s-t)^\alpha}\, ds. \qquad (2.2)$$

In addition, the right Riemann-Liouville differentiation operator is written in the form

$${}^{RL}_{t}\mathcal{D}^{\alpha}_{B}\varphi(t) = -\frac{d}{dt}\int_{t}^{B}\frac{\varphi(s)}{(s-t)^{\alpha}}\,ds. \qquad (2.3)$$

On the other hand, the fractional integral is given by

$${}_{A}\mathcal{I}^{\alpha}_{t}\varphi(t) = \frac{1}{\Gamma(\alpha)}\int_{A}^{t}(t-s)^{\alpha-1}\varphi(s)\,ds. \qquad (2.4)$$

Lastly, we mention a useful relation between right Caputo and Riemann-Liouville differentiation operators

$${}^{RL}_{t}\mathcal{D}^{\alpha}_{B}\varphi(t) = {}^{C}_{t}\mathcal{D}^{\alpha}_{B}\varphi(t) + \varphi(B)\frac{(B-t)^{-\alpha}}{\Gamma(1-\alpha)}. \qquad (2.5)$$

In this study, we consider a FOCP governed by a generalized cancer-obesity model. The model describes the interaction between the state variables, namely tumor cells $T(t)$, immune cells $I(t)$, fat cells $F(t)$, chemotherapeutic and immunotherapeutic drug concentrations $D_1(t)$ and $D_2(t)$, respectively. The control variables $u_1(t)$ and $u_2(t)$ denote the doses of the chemotherapeutic and immunotherapeutic drugs, respectively. The aim behind the OCP is to minimize the value of the cost functional $J(u_1, u_2)$, which is equivalent to minimizing the tumor population and the drug administered over a finite time interval, where the interaction of tumor cells, immune cells, fat cells and drug concentrations are governed by a time fractional coupled system of differential equations.

We shortly denote

$$T := T(t),\ I := I(t),\ F := F(t),\ D_1 := D_1(t),\ D_2 := D_2(t),$$

and propose the following FOCP

$$\min_{u=(u_1,u_2)\in U_{ad}} J(u_1, u_2) = \int_0^{t_f} (T + \omega_1 u_1^2 + \omega_2 u_2^2)\,dt \qquad (2.6)$$

subject to

$$\begin{cases} {}^{C}_{0}\mathcal{D}^{\alpha}_{t}T &= r^{\alpha}T(1-p^{\alpha}T) - \xi_1^{\alpha}TI + c_1^{\alpha}TF - q_1^{\alpha}D_1T, \\ {}^{C}_{0}\mathcal{D}^{\alpha}_{t}I &= s^{\alpha} + \frac{\rho^{\alpha}T^2 I}{h^{\alpha}+T^2+F^2} + \frac{\beta^{\alpha}D_2 I}{g^{\alpha}+D_2} - \xi_2^{\alpha}TI - \mu^{\alpha}I - q_2^{\alpha}D_1 I, \\ {}^{C}_{0}\mathcal{D}^{\alpha}_{t}F &= d^{\alpha}F(1-\epsilon^{\alpha}F) - c_2^{\alpha}FT - q_3^{\alpha}D_1 F, \\ {}^{C}_{0}\mathcal{D}^{\alpha}_{t}D_1 &= u_1 - \gamma_1^{\alpha}D_1, \\ {}^{C}_{0}\mathcal{D}^{\alpha}_{t}D_2 &= u_2 - \gamma_2^{\alpha}D_2, \\ T(0) &= T_0, \ I(0) = I_0, \ F(0) = F_0, \ D_1(0) = D_{10}, \ D_2(0) = D_{20}, \end{cases} \quad (2.7)$$

where the admissible space of controls is given by

$$U_{ad} = \{u = (u_1, u_2) \mid u_1, u_2 \text{ are measurable with } 0 \leq u_1, u_2 \leq 1, \ t \in [0, t_f]\}.$$

Therefore, the aim is to find the optimal control $u^* = (u_1^*, u_2^*) \in U_{ad}$ such that $J(u_1^*, u_2^*) = \min_{u=(u_1,u_2)\in U_{ad}} J(u_1, u_2)$ holds.

Now, we proceed with the positivity of the solution of the state equation (2.7). Then, we obtain the optimality system for the OCP (2.6)-(2.7).

## 2.1 Existence of positive solutions of the cancer-obesity model

We denote $\mathbb{R}_+^5 = \{x \in \mathbb{R}^5 \mid x \geq 0\}$ and $x(t) = (T, I, F, D_1, D_2)^T$.

**Theorem 2.1.** *The solution of the FDE (2.7) is unique and it remains in $\mathbb{R}_+^5$.*

*Proof.* From [19, Theorem 3.1, Remark 3.2], existence of the unique solution to the FDE (2.7) is proven on $(0, \infty)$. Then, we show that the non-negative orthant $\mathbb{R}_+^5$ is a positively invariant region. Since,

$$\begin{gathered} {}^{C}_{0}\mathcal{D}^{\alpha}_{t}T|_{T=0} = 0, \quad {}^{C}_{0}\mathcal{D}^{\alpha}_{t}I|_{I=0} = s^{\alpha} \geq 0, \quad {}^{C}_{0}\mathcal{D}^{\alpha}_{t}F|_{F=0} = 0, \\ {}^{C}_{0}\mathcal{D}^{\alpha}_{t}D_1|_{D_1=0} = u_1 \geq 0, \quad {}^{C}_{0}\mathcal{D}^{\alpha}_{t}D_2|_{D_2=0} = u_2 \geq 0, \end{gathered}$$

on each line bounding the non-negative orthant, the vector field points into $\mathbb{R}_+^5$. The solution will remain in $\mathbb{R}_+^5$ [23]. □

## 2.2 Optimality system

We derive the necessary optimality conditions for the OCP (2.6-2.7).

**Theorem 2.2.** *Given a pair of optimal controls $u^* = (u_1^*, u_2^*)$ and the state solutions $(T^*, I^*, F^*, D_1^*, D_2^*)$ corresponding to (2.7) that minimize objective functional (2.6), there* *exist adjoint variables* $(\lambda_1, \lambda_2, \lambda_3, \lambda_4, \lambda_5)$ *satisfying*

$$\begin{cases} {}^C_t\mathcal{D}^\alpha_{t_f}\lambda_1 &= (r^\alpha - 2r^\alpha p^\alpha T - \xi_1^\alpha I + c_1^\alpha F - q_1^\alpha D_1)\lambda_1 \\ &\quad + \left(\frac{2\rho^\alpha TI(h^\alpha+F^2)}{(h^\alpha+T^2+F^2)^2} - \xi_2^\alpha I\right)\lambda_2 + 1, \\ {}^C_t\mathcal{D}^\alpha_{t_f}\lambda_2 &= -\xi_1^\alpha T\lambda_1 + \left(\frac{\rho^\alpha T^2}{h^\alpha+T^2+F^2} + \frac{\beta^\alpha D_2}{g^\alpha+D_2} - \xi_2^\alpha T - \mu^\alpha - q_2^\alpha D_1\right)\lambda_2, \\ {}^C_t\mathcal{D}^\alpha_{t_f}\lambda_3 &= c_1^\alpha T\lambda_1 - \frac{2\rho^\alpha T^2 IF}{(h^\alpha+T^2+F^2)^2}\lambda_2 + (d^\alpha - 2d^\alpha\epsilon^\alpha F - c_2^\alpha T - q_3^\alpha D_1)\lambda_3, \\ {}^C_t\mathcal{D}^\alpha_{t_f}\lambda_4 &= -q_1^\alpha T\lambda_1 - q_2^\alpha I\lambda_2 - q_3^\alpha F\lambda_3 - \gamma_1^\alpha \lambda_4, \\ {}^C_t\mathcal{D}^\alpha_{t_f}\lambda_5 &= \frac{\beta^\alpha g^\alpha I}{(g^\alpha+D_2)^2}\lambda_2 - \gamma_2^\alpha \lambda_5, \end{cases} \quad (2.8)$$

*with transversality conditions*

$$\lambda_1(t_f) = 0, \ \lambda_2(t_f) = 0, \ \lambda_3(t_f) = 0, \ \lambda_4(t_f) = 0, \ \lambda_5(t_f) = 0. \quad (2.9)$$

*Moreover, the pair of optimal controls* $u^* = (u_1^*, u_2^*)$ *is represented by*

$$\begin{cases} u_1^* &= \min\left(\max\left(-\frac{\lambda_4}{2\omega_1}, 0\right), 1\right), \\ u_2^* &= \min\left(\max\left(-\frac{\lambda_5}{2\omega_2}, 0\right), 1\right). \end{cases} \quad (2.10)$$

*Proof.* Following the proof in [36, Theorem 5.1], we note that existence of a pair of optimal controls $u^* = (u_1^*, u_2^*)$ and the associated state solution $(T^*, I^*, F^*, D_1^*, D_2^*)$ is obtained due to the convexity of $J(u_1, u_2)$ with respect to controls and the constraint, which satisfies Lipschitz property with respect to state variables.

The optimality system can be derived by constructing the Lagrangian as

$$\begin{aligned} &\mathcal{L}(T, I, F, D_1, D_2, u_1, u_2, \lambda_1, \lambda_2, \lambda_3, \lambda_4, \lambda_5) \\ &= \int_0^{t_f} \Big((T + \omega_1 u_1^2 + \omega_2 u_2^2) \\ &\quad - \lambda_1^T\left({}^C_0\mathcal{D}^\alpha_t T - r^\alpha T(1-p^\alpha T) + \xi_1^\alpha TI - c_1^\alpha TF + q_1^\alpha D_1 T\right) \\ &\quad - \lambda_2^T\left({}^C_0\mathcal{D}^\alpha_t I - s^\alpha - \frac{\rho^\alpha T^2 I}{h^\alpha+T^2+F^2} - \frac{\beta^\alpha D_2 I}{g^\alpha+D_2} + \xi_2^\alpha TI + \mu^\alpha I + q_2^\alpha D_1 I\right) \\ &\quad - \lambda_3^T\left({}^C_0\mathcal{D}^\alpha_t F - d^\alpha F(1-\epsilon^\alpha F) + c_2^\alpha FT + q_3^\alpha D_1 F\right) \\ &\quad - \lambda_4^T\left({}^C_0\mathcal{D}^\alpha_t D_1 - u_1 + \gamma_1^\alpha D_1\right) - \lambda_5^T\left({}^C_0\mathcal{D}^\alpha_t D_2 - u_2 + \gamma_2^\alpha D_2\right)\Big) dt \\ &\quad + \xi_1 u_1 + \xi_2(1-u_1) + \mu_1 u_2 + \mu_2(1-u_2) \\ &\quad - \lambda_1(0)(T(0) - T_0) - \lambda_2(0)(I(0) - I_0) - \lambda_3(0)(F(0) - F_0) \\ &\quad - \lambda_4(0)(D_1(0) - D_{10}) - \lambda_5(0)(D_2(0) - D_{20}), \end{aligned} \quad (2.11)$$

where $\lambda_i(t)$'s are the adjoint or co-state variables, $\xi_1 \geq 0, \xi_2 \geq 2, \mu_1 \geq 0, \mu_2 \geq 0$ are penalty multipliers satisfying

$$\xi_1 u_1 = 0, \quad \xi_2(1-u_1) = 0, \quad \mu_1 u_2 = 0, \quad \mu_2(1-u_2) = 0,$$

at the optimal $u^* = (u_1^*, u_2^*)$. In particular, we present the derivation of the fractional derivative of the adjoint equation. We apply the method of integration by parts, for example, to the term $\lambda_1^T(t)\left({}_0^C\mathcal{D}_t^\alpha T(t)\right)$ in (2.11) following the study [1, Sec. 2] as

$$\int_0^{t_f} \lambda_1^T(t)\left({}_0^C\mathcal{D}_t^\alpha T(t)\right) dt = \int_0^{t_f} \left({}_t^R\mathcal{D}_{t_f}^\alpha \lambda_1(t)\right)^T T(t)\, dt, \qquad (2.12)$$

with the condition $\lambda_1(t_f) = 0$. Then, we proceed as

$$\int_0^{t_f} \lambda_1^T(t)\left({}_t^R\mathcal{D}_{t_f}^\alpha T(t)\right) dt = \int_0^{t_f} \left({}_t^C\mathcal{D}_{t_f}^\alpha \lambda_1(t)\right)^T T(t)\, dt$$

$$\underset{(2.5)}{=} \int_0^{t_f} \left({}_t^C\mathcal{D}_{t_f}^\alpha \lambda_1(t)\right)^T T(t)\, dt + \overbrace{\frac{\lambda_1(t_f)(t_f - t)^{-\alpha}}{\Gamma(1-\alpha)}}^{=0}$$

$$= \int_0^{t_f} \left({}_t^C\mathcal{D}_{t_f}^\alpha \lambda_1(t)\right)^T T(t)\, dt. \qquad (2.13)$$

Then, we substitute (2.13) into (2.11) and differentiate the resulting equation with respect to $(T, I, F, D_1, D_2)$ and $(u_1, u_2)$ to obtain the adjoint equation (2.8) with final conditions (2.9) and the gradient equation (2.10), respectively. $\square$

## 3   Equilibrium points and stability analysis

In this section, we will examine the existence of the equilibrium points of the system (2.7), namely, *Tumor Free Equilibrium* and *Coexisting Equilibrium*, and determine the conditions under which they are stable. Firstly, we fix the doses of chemotherapeutic and immunotherapeutic drugs as $u_1(t) = u_1$ and $u_2(t) = u_2$, respectively. Then, we solve the following system of FDEs for $T(t), I(t), F(t), D_1(t)$ and $D_2(t)$:

$$ {}_0^C\mathcal{D}_t^\alpha T = 0, \quad {}_0^C\mathcal{D}_t^\alpha I = 0, \quad {}_0^C\mathcal{D}_t^\alpha F = 0, \quad {}_0^C\mathcal{D}_t^\alpha D_1 = 0, \quad {}_0^C\mathcal{D}_t^\alpha D_2 = 0. \qquad (3.1)$$

## 3.1 Tumor free equilibrium

We consider the case that no tumor cells exist, that is, $\hat{T} = 0$. We solve the system (3.1) and find the tumor-free equilibrium point $\hat{E} = (0, \hat{I}, \hat{F}, \hat{D}_1, \hat{D}_2)$ where

$$\hat{I} = \frac{s^\alpha(g^\alpha + \hat{D}_2)}{(\mu^\alpha + q_2^\alpha \hat{D}_1)(g^\alpha + \hat{D}_2) - \beta^\alpha \hat{D}_2},$$

$$\hat{F} = \frac{1}{\epsilon^\alpha} - \frac{q_3^\alpha}{d^\alpha \epsilon^\alpha}\hat{D}_1, \qquad (3.2)$$

$$\hat{D}_1 = \frac{u_1}{\gamma_1^\alpha}, \quad \hat{D}_2 = \frac{u_2}{\gamma_2^\alpha}.$$

The equilibrium point $\hat{E}$ exists if $\hat{I} > 0$, $\hat{F} > 0$, $\hat{D}_1 > 0$ and $\hat{D}_2 > 0$. We note that $\hat{D}_1 > 0$ and $\hat{D}_2 > 0$ are automatically satisfied with $u_1 > 0$ and $u_2 > 0$ according to (3.2). Now, we must assure that $\hat{I} > 0$ and $\hat{F} > 0$ hold. Then, we obtain the following inequalities

$$\hat{I} = \frac{s^\alpha(g^\alpha + \hat{D}_2)}{(\mu^\alpha + q_2^\alpha \hat{D}_1)(g^\alpha + \hat{D}_2) - \beta^\alpha \hat{D}_2}$$

$$= \frac{s^\alpha \gamma_1^\alpha (\gamma_2^\alpha g^\alpha + u_2)}{\mu^\alpha g^\alpha \gamma_1^\alpha \gamma_2^\alpha + q_2^\alpha g^\alpha u_1 \gamma_2^\alpha + q_2^\alpha u_1 u_2 + u_2 \gamma_1^\alpha (\mu^\alpha - \beta^\alpha)} > 0,$$

and

$$\hat{F} = \frac{1}{\epsilon^\alpha} - \frac{q_3^\alpha}{d^\alpha \epsilon^\alpha}\hat{D}_1 = \frac{1}{\epsilon^\alpha} - \frac{q_3^\alpha}{d^\alpha \epsilon^\alpha}\frac{u_1}{\gamma_1^\alpha} > 0.$$

Therefore, we find the following conditions

$$\mu^\alpha g^\alpha \gamma_1^\alpha \gamma_2^\alpha + q_2^\alpha g^\alpha u_1 \gamma_2^\alpha + q_2^\alpha u_1 u_2 + u_2 \gamma_1^\alpha \mu^\alpha > u_2 \gamma_1^\alpha \beta^\alpha, \qquad d^\alpha \gamma_1^\alpha > q_3^\alpha u_1,$$

so that the equilibrium point $\hat{E}$ exists.

Now, we discuss the stability of $\hat{E}$ by investigating the signs of the eigenvalues of the Jacobian associated to the model (2.7).

**Theorem 3.1.** *The equilibrium point $\hat{E} = (0, \hat{I}, \hat{F}, \hat{D}_1, \hat{D}_2)$ of the system (2.7) exists under the condition that*

$$\mu^\alpha g^\alpha \gamma_1^\alpha \gamma_2^\alpha + q_2^\alpha g^\alpha u_1 \gamma_2^\alpha + q_2^\alpha u_1 u_2 + u_2 \gamma_1^\alpha \mu^\alpha > u_2 \gamma_1^\alpha \beta^\alpha, \qquad d^\alpha \gamma_1^\alpha > q_3^\alpha u_1. \qquad (3.3)$$

Moreover, $\hat{E} = (0, \hat{I}, \hat{F}, \hat{D}_1, \hat{D}_2)$ is locally asymptotically stable if the following conditions

$$
\begin{aligned}
& q_3^\alpha u_1 < d^\alpha \gamma_1^\alpha, \\
& \beta^\alpha \gamma_1^\alpha u_2 < (\gamma_1^\alpha \mu^\alpha + q_2^\alpha u_1)(\gamma_2^\alpha g^\alpha + u_2), \\
& I^* > \frac{r^\alpha}{\xi_1^\alpha} + \frac{c_1^\alpha}{\xi_1^\alpha \epsilon^\alpha} - \left( \frac{c_1^\alpha q_3^\alpha}{\xi_1^\alpha d^\alpha \epsilon^\alpha} + \frac{q_1^\alpha}{\xi_1^\alpha \gamma_1^\alpha} \right) u_1,
\end{aligned}
\tag{3.4}
$$

*together with* (3.3) *hold.*

*Proof.* The Jacobian of the system (2.7) evaluated at $\hat{E} = (0, \hat{I}, \hat{F}, \hat{D}_1, \hat{D}_2)$ is given by

$$
J(E^*) = \begin{pmatrix}
\omega_1 & 0 & 0 & 0 & 0 \\
\omega_2 & \omega_3 & 0 & \omega_4 & \omega_5 \\
\omega_6 & 0 & \omega_7 & \omega_8 & 0 \\
0 & 0 & 0 & \omega_9 & 0 \\
0 & 0 & 0 & 0 & \omega_{10}
\end{pmatrix},
$$

where

$$
\begin{aligned}
& \omega_1 = r^\alpha - \xi_1^\alpha \hat{I} + c_1^\alpha \hat{F} - q_1^\alpha \hat{D}_1, \quad \omega_2 = -\xi_2^\alpha \hat{I}, \\
& \omega_3 = -\mu^\alpha - q_2^\alpha \hat{D}_1 + \beta^\alpha \frac{\hat{D}_2}{g^\alpha + \hat{D}_2}, \quad \omega_4 = -q_2^\alpha \hat{I}, \quad \omega_5 = \frac{g^\alpha \beta^\alpha \hat{I}}{(g^\alpha + \hat{D}_2)^2}, \\
& \omega_6 = -c_2^\alpha \hat{F}, \quad \omega_7 = d^\alpha - 2d^\alpha \epsilon^\alpha \hat{F} - q_3^\alpha \hat{D}_1, \\
& \omega_8 = -q_3^\alpha \hat{F}, \quad \omega_9 = -\gamma_1^\alpha, \quad \omega_{10} = -\gamma_2^\alpha.
\end{aligned}
$$

The eigenvalues of $J(\hat{E})$ are computed as

$$
\begin{aligned}
\lambda_1 &= \omega_9 = -\gamma_1^\alpha, \\
\lambda_2 &= \omega_{10} = -\gamma_2^\alpha, \\
\lambda_3 &= \omega_3 = \frac{\beta^\alpha u_2}{\gamma_2^\alpha g^\alpha + u_2} - \left( \mu^\alpha + q_2^\alpha \frac{u_1}{\gamma_1^\alpha} \right), \\
\lambda_4 &= \omega_7 = \frac{q_3^\alpha u_1}{\gamma_1^\alpha} - d^\alpha, \\
\lambda_5 &= \omega_1 = r^\alpha - \xi_1^\alpha \hat{I} + c_1^\alpha \hat{F} - q_1^\alpha \hat{D}_1 \\
&= -\xi_1^\alpha \hat{I} + r^\alpha + \frac{c_1^\alpha}{\epsilon^\alpha} - \left( \frac{c_1^\alpha q_3^\alpha}{d^\alpha \epsilon^\alpha} + \frac{q_1^\alpha}{\gamma_1^\alpha} \right) u_1.
\end{aligned}
$$

Tumor-free equilibrium point $\hat{E}$ is locally asymptotically stable if all eigenvalues $\lambda_i$ of $J(\hat{E})$ are negative. The requirement $\lambda_1 < 0$ and $\lambda_2 < 0$ are automatically satisfied. For $\lambda_3, \lambda_4, \lambda_5$ to be negative, we reach the conditions in the Eqn. (3.4). $\square$

## 3.2 Co-existing equilibrium point

In this case, tumor cells coexist, i.e. $\overline{T} \neq 0$. By keeping this in mind, we solve (3.1) for $\overline{T} > 0$, $\overline{I} > 0$, $\overline{F} > 0$, $\overline{D}_1 > 0$, $\overline{D}_2 > 0$ to reach the co-existing equilibrium point $\overline{E} = (\overline{T}, \overline{I}, \overline{F}, \overline{D}_1, \overline{D}_2)$.

We find that

$$\overline{D}_1 = \frac{u_1}{\gamma_1^\alpha}, \quad \overline{D}_2 = \frac{u_2}{\gamma_2^\alpha},$$
$$\overline{F} = \frac{1}{\epsilon^\alpha d^\alpha}(d^\alpha - q_3^\alpha \overline{D}_1 - c_2^\alpha \overline{T}) = k_1 + k_2 \overline{T}, \qquad (3.5)$$
$$\overline{T} = \frac{1}{p^\alpha r^\alpha}(r^\alpha - \xi_1^\alpha \overline{I} - q_1^\alpha \overline{D}_1 + c_1^\alpha \overline{F}) = k_4 + k_5 \overline{I},$$

where

$$k_1 = \frac{1}{\epsilon^\alpha} - \frac{q_3^\alpha}{\epsilon^\alpha d^\alpha}\overline{D}_1, \quad k_2 = -\frac{c_2^\alpha}{\epsilon^\alpha d^\alpha}, \quad k_3 = \frac{1}{p^\alpha} - \frac{q_1^\alpha}{p^\alpha r^\alpha}\overline{D}_1,$$
$$k_4 = \frac{k_3^\alpha p^\alpha r^\alpha + k_1 c_1^\alpha}{p^\alpha r^\alpha - c_1^\alpha k_2} = \frac{\epsilon^\alpha d^\alpha (k_3^\alpha p^\alpha r^\alpha + k_1 c_1^\alpha)}{\epsilon^\alpha d^\alpha p^\alpha r^\alpha + c_1^\alpha c_2^\alpha},$$
$$k_5 = -\frac{\xi_1^\alpha}{p^\alpha r^\alpha - c_1^\alpha k_2} = -\frac{\epsilon^\alpha d^\alpha \xi_1^\alpha}{\epsilon^\alpha d^\alpha p^\alpha r^\alpha + c_1^\alpha c_2^\alpha}.$$

For simplification purposes, we put

$$\overline{F} = k_1 + k_2(k_4 + k_5\overline{I}) = k_6 + k_7\overline{I},$$
$$\overline{I} = \frac{s^\alpha(h^\alpha + \overline{T}^2 + \overline{F}^2)(g^\alpha + \overline{D}_2)}{(\mu^\alpha + \alpha_2^\alpha \overline{T} + q_2^\alpha \overline{D}_1)(h^\alpha + \overline{T}^2 + \overline{F}^2)(g^\alpha + \overline{D}_2) - \overline{T}^2(\rho^\alpha(g^\alpha + \overline{D}_2) + \beta^\alpha \overline{D}_2) - \beta^\alpha \overline{D}_2(\overline{F}^2 + h^\alpha)},$$

where $k_6 = k_1 + k_2 k_4$, $k_7 = k_2 k_5$. After arranging the terms, we obtain the following polynomial in $\overline{I}$

$$m_4 \overline{I}^4 + m_3 \overline{I}^3 + m_2 \overline{I}^2 + m_1 \overline{I} + m_0 = 0, \qquad (3.6)$$

where

$$m_4 = \theta_1 k_5 (k_5^2 + k_7^2),$$
$$m_3 = \theta_1(3k_4 k_5^2 + k_4 k_7^2 + 2k_5 k_6 k_7) + \theta_2 k_7^2 + \theta_3 k_5^2,$$
$$m_2 = \theta_1(3k_4^2 k_5 + 2k_4 k_6 k_7 + k_5 k_6^2) + \theta_2(2k_6 k_7) + \theta_3(2k_4 k_5) + \theta_4 k_5 + \theta_5(k_5^2 + k_7^2),$$
$$m_1 = \theta_1(k_4^3 + k_4 k_6^2) + \theta_2 k_6^2 + \theta_3 k_5^2 + \theta_4 k_4 + \theta_5(2k_4 k_5 + 2k_6 k_7) + \theta_6,$$
$$m_0 = \theta_5(k_4^2 + k_6^2) + \theta_7,$$
$$\theta_1 = (g^\alpha + \overline{D}_2)\xi_2^\alpha,$$
$$\theta_2 = (g^\alpha + \overline{D}_2)(\mu^\alpha + q_2^\alpha \overline{D}_1) - \beta^\alpha \overline{D}_2,$$
$$\theta_3 = (g^\alpha + \overline{D}_2)(\mu^\alpha + q_2^\alpha \overline{D}_1) - (\rho^\alpha g^\alpha + \rho^\alpha \overline{D}_2 + \beta^\alpha \overline{D}_2),$$
$$\theta_4 = (g^\alpha + \overline{D}_2)\xi_2^\alpha h^\alpha,$$
$$\theta_5 = -(g^\alpha + \overline{D}_2)s^\alpha,$$
$$\theta_6 = (g^\alpha + \overline{D}_2)(h^\alpha \mu^\alpha + q_2^\alpha h^\alpha \overline{D}_1) - \beta^\alpha h^\alpha \overline{D}_2,$$
$$\theta_7 = -(g^\alpha + \overline{D}_2)s^\alpha h^\alpha.$$

We observe that

$$m_4 = \theta_1 k_5(k_5^2 + k_7^2) < 0,$$
$$m_0 = \theta_5(k_4^2 + k_6^2) + \theta_7 = -(g^\alpha + \overline{D}_2)s^\alpha(k_4^2 + k_6^2 + h^\alpha) < 0.$$

We obtain following conditions on $m_1, m_2, m_3$ through the use of Descartes's rule so that the Eqn. (3.6) may have non-trivial positive roots:

(1) if $m_1 < 0$, $m_2 < 0$, $m_3 < 0$, then there is no change of sign, so no positive roots of Eqn. (3.6) exists,

(2) if $m_1 > 0$, $m_2 < 0$, $m_3 < 0$, two or no positive roots of Eqn. (3.6) exist,

(3) if $m_1 < 0$, $m_2 > 0$, $m_3 < 0$, two or no positive roots of Eqn. (3.6) exist,

(4) if $m_1 > 0$, $m_2 > 0$, $m_3 < 0$, two or no positive roots of Eqn. (3.6) exist,

(5) if $m_1 < 0$, $m_2 < 0$, $m_3 > 0$, two or no positive roots of Eqn. (3.6) exist,

(6) if $m_1 < 0$, $m_2 > 0$, $m_3 > 0$, two or no positive roots of Eqn. (3.6) exist,

(7) if $m_1 > 0$, $m_2 > 0$, $m_3 > 0$, two or no positive roots of Eqn. (3.6) exist,

(8) if $m_1 > 0$, $m_2 < 0$, $m_3 > 0$, four or two or no positive roots of Eqn. (3.6) exist.

Hence, $\overline{I}$ is not trivially positive if any one of the above seven conditions $(2) - (8)$ is satisfied. We note that $\overline{D}_1 > 0$ and $\overline{D}_2 > 0$ hold due to the non-negative dosage of drugs $u_1 > 0$ and $u_2 > 0$. Now, we must find the conditions so that $\overline{T} > 0$ and $\overline{F} > 0$ hold. According to the equations in (3.5), we find that

$$\overline{F} > 0 \text{ if } k_1 + k_2\overline{T} > 0, \text{ i.e., } \overline{T} < -\frac{k_1}{k_2} = \frac{d^\alpha \gamma_1^\alpha - q_3^\alpha u_1}{\gamma_1^\alpha c_2^\alpha},$$

$$\overline{T} > 0 \text{ if } k_4 + k_5\overline{I} > 0, \text{ i.e., } \overline{I} < -\frac{k_4}{k_5} = \frac{k_3 p^\alpha r^\alpha + k_1 c_1^\alpha}{\xi_1^\alpha}.$$

We summarize the conditions on existence of the equilibrium point $\overline{E} = (\overline{T}, \overline{I}, \overline{F}, \overline{D}_1, \overline{D}_2)$ in the following theorem:

**Theorem 3.2.** *The equilibrium point $\overline{E} = (\overline{T}, \overline{I}, \overline{F}, \overline{D}_1, \overline{D}_2)$ of the system (2.7) exists under the condition that*

$$\overline{T} < \frac{d^\alpha \gamma_1^\alpha - q_3^\alpha u_1}{\gamma_1^\alpha c_2^\alpha}, \qquad \overline{I} < \frac{k_3 p^\alpha r^\alpha + k_1 c_1^\alpha}{\xi_1^\alpha}, \tag{3.7}$$

*and any one of the seven conditions $(2)-(8)$ above is satisfied, i.e., at least one of $m_1, m_2, m_3$ is positive.*

We proceed with investigation of the stability of the co-existing equilibrium point $\overline{E} = (\overline{T}, \overline{I}, \overline{F}, \overline{D}_1, \overline{D}_2)$. We note that if the eigenvalues, namely $\lambda_i$'s, of the Jacobian matrix evaluated at $\overline{E}$ satisfy the condition

$$|arg(\lambda_i)| > \frac{\alpha \pi}{2}, \quad (i = 1, 2, 3, 4, 5), \tag{3.8}$$

then the system is asymptotically stable at $\overline{E}$ [21, 2]. In detail, Jacobian matrix of the system (2.7) evaluated at $\overline{E} = (\overline{T}, \overline{I}, \overline{F}, \overline{D}_1, \overline{D}_2)$ is given by

$$J(\overline{E}) = \begin{pmatrix} \omega_1 & \omega_2 & \omega_3 & \omega_4 & 0 \\ \omega_5 & \omega_6 & \omega_7 & \omega_8 & \omega_9 \\ \omega_{10} & 0 & \omega_{11} & \omega_{12} & 0 \\ 0 & 0 & 0 & \omega_{13} & 0 \\ 0 & 0 & 0 & 0 & \omega_{14} \end{pmatrix},$$

where

$$\omega_1 = r^\alpha - 2r^\alpha p^\alpha \overline{T} - \xi_1^\alpha \overline{I} + c_1^\alpha \overline{F} - q_1^\alpha \overline{D}_1, \quad \omega_2 = -\xi_1^\alpha \overline{T},$$

$$\omega_3 = c_1^\alpha \overline{T}, \quad \omega_4 = -q_1^\alpha \overline{T}, \quad \omega_5 = \frac{2\rho^\alpha \overline{TI}(h^\alpha + \overline{F}^2)}{(h^\alpha + \overline{T}^2 + \overline{F}^2)^2} - \xi_2 \overline{I},$$

$$\omega_6 = \frac{\rho^\alpha \overline{T}^2}{h^\alpha + \overline{T}^2 + \overline{F}^2} + \frac{\beta^\alpha \overline{D}_2}{g^\alpha + \overline{D}_2} - \xi_2^\alpha \overline{T} - \mu^\alpha - q_2^\alpha \overline{D}_1,$$

$$\omega_7 = -\frac{2\rho^\alpha \overline{T}^2 \overline{IF}}{(h^\alpha + \overline{T}^2 + \overline{F}^2)^2}, \quad \omega_8 = -q_2^\alpha \overline{I}, \quad \omega_9 = \frac{\beta^\alpha g^\alpha \overline{I}}{(g^\alpha + \overline{D}_2)^2},$$

$$\omega_{10} = -c_2^\alpha \overline{F}, \quad \omega_{11} = d^\alpha - 2d^\alpha \epsilon^\alpha \overline{F} - c_2^\alpha \overline{T} - q_3^\alpha \overline{D}_1,$$

$$\omega_{12} = -q_3^\alpha \overline{F}, \quad \omega_{13} = -\gamma_1^\alpha, \quad \omega_{14} = -\gamma_2^\alpha.$$

The eigenvalues of $J(\overline{E})$ are the roots of the following characteristic polynomial

$$p(\lambda) = -(\lambda - \omega_{13})(\lambda - \omega_{14}) \underbrace{(\lambda^3 + c_1 \lambda^2 + c_2 \lambda + c_3)}_{:=q(\lambda)}, \tag{3.9}$$

where

$$\begin{aligned}
c_1 &= (\omega_1 + \omega_6 + \omega_{11}), \\
c_2 &= (-\omega_1 \omega_6 + \omega_2 \omega_5 - \omega_1 \omega_{11} + \omega_3 \omega_{10} - \omega_6 \omega_{11}), \\
c_3 &= (\omega_1 \omega_6 \omega_{11} - \omega_2 \omega_5 \omega_{11} + \omega_2 \omega_7 \omega_{10} - \omega_3 \omega_6 \omega_{10}).
\end{aligned} \tag{3.10}$$

We can directly deduce that two of the eigenvalues satisfy $\lambda_1 = \omega_{13} = -\gamma_1^\alpha < 0$ and $\lambda_2 = \omega_{14} = -\gamma_2^\alpha < 0$. Therefore, we proceed with the roots of the polynomial $q(\lambda)$. We note that the discriminant of $q(\lambda)$ is given by

$$D(q) = 18 c_1 c_2 c_3 + (c_1 c_2)^2 - 4 c_3 c_1^3 - 4 c_2^3 - 27 c_3^2.$$

Then, we list the following conditions for the model (2.7) so that the condition (3.8) holds according to Routh–Hurwitz criteria [2]:

**Corollary 3.3.** *Assume that the conditions of Theorem 3.2 hold so that the co-existing equilibrium point exists. Then, the equilibrium point $\overline{E} = (\overline{T}, \overline{I}, \overline{F}, \overline{D}_1, \overline{D}_2)$ is locally asymptotically stable if one of the following conditions holds for polynomial $q(\lambda)$ which is given as in (3.9) and coefficients $c_1, c_2, c_3$ which are given as in (3.10).*

(i) *If $D(q) > 0$, then the necessary and sufficient condition for the equilibrium point $\overline{E}$ to be locally asymptotically stable is $c_1 > 0$, $c_3 > 0$, $c_1 c_2 > c_3$.*

(ii) If $D(q) < 0$, $c_1 \geq 0$, $c_2 \geq 0$, $c_3 > 0$, then $\overline{E}$ is locally asymptotically stable for $\alpha < 2/3$.

(iii) If $D(q) < 0$, $c_1 > 0$, $c_2 > 0$, $c_1 c_2 = c_3$, then $\overline{E}$ is locally asymptotically stable.

# 4 Discretization technique

Let $0 = t_0 < t_1 < \ldots < t_N = t_f$ be a subdivision of $I = (0, t_f]$ with constant time step $\Delta t = T/N$. We denote the approximate value of $\varphi(t)$ at $t = t_j$ as $\varphi_j$.

We present the derivation of the discrete (left/right) Caputo fractional derivative.

## 4.1 Discrete state equation

L1-method for (left) Caputo fractional derivative has presented in the paper [20, Sec. 3]. We follow the same idea to obtain

$$\begin{aligned}
{}_0^C\mathbb{D}_t^\alpha \varphi(t)|_{t=t_k} &= \frac{1}{\Gamma(1-\alpha)} \sum_{j=1}^{k} \frac{\varphi_j - \varphi_{j-1}}{\Delta t} \int_{t_{j-1}}^{t_j} (t_k - s)^{-\alpha}\, ds \\
&= B_0 \sum_{j=1}^{k} (\varphi_j - \varphi_{j-1})\left(\delta_{j,k}^C\right),
\end{aligned} \quad (4.1)$$

where $B_0 = \frac{-\Delta t^{-\alpha}}{\Gamma(2-\alpha)}$, $\delta_{j,k}^C = ((k-j)^{1-\alpha} - (k-j+1)^{1-\alpha})$.

We apply the scheme (4.1) to the state equation (2.7). Then, the nonlinear state equation is linearized by Newton's method for $1 \leq k \leq N$ and the solution is obtained by solving the following system iteratively:

$$\begin{pmatrix}
\delta_{k,k}^C B_0 - J_{11} & -J_{12} & -J_{13} & -J_{14} & 0 \\
-J_{21} & \delta_{k,k}^C B_0 - J_{22} & -J_{23} & -J_{24} & -J_{25} \\
-J_{31} & 0 & \delta_{k,k}^C B_0 - J_{33} & -J_{34} & 0 \\
0 & 0 & 0 & \delta_{k,k}^C B_0 - J_{44} & 0 \\
0 & 0 & 0 & 0 & \delta_{k,k}^C B_0 - J_{55}
\end{pmatrix}
\begin{pmatrix} \delta T \\ \delta I \\ \delta F \\ \delta D_1 \\ \delta D_2 \end{pmatrix}
=
\begin{pmatrix} R_1 \\ R_2 \\ R_3 \\ R_4 \\ R_5 \end{pmatrix}, \quad (4.2)$$

where

$$\begin{cases}
R_1 = {}_0^C\mathbb{D}_t^\alpha T(t)|_{t=t_k} - (r^\alpha T_k(1 - p^\alpha T_k) - \xi_1^\alpha T_k I_k + c_1^\alpha T_k F_k - q_1^\alpha (D_1)_k T_k), \\
R_2 = {}_0^C\mathbb{D}_t^\alpha I(t)|_{t=t_k} - (s^\alpha + \frac{\rho^\alpha T_k^2 I_k}{h^\alpha + T_k^2 + F_k^2} + \frac{\beta^\alpha (D_2)_k I_k}{g^\alpha + (D_2)_k} - \xi_2^\alpha T_k I_k - \mu^\alpha I_k - q_2^\alpha (D_1)_k I_k), \\
R_3 = {}_0^C\mathbb{D}_t^\alpha F(t)|_{t=t_k} - (d^\alpha F_k(1 - \epsilon^\alpha F_k) - c_2^\alpha F_k T_k - q_3^\alpha (D_1)_k F_k), \\
R_4 = {}_0^C\mathbb{D}_t^\alpha D_1(t)|_{t=t_k} - ((u_1)_k - \gamma_1^\alpha (D_1)_k), \\
R_5 = {}_0^C\mathbb{D}_t^\alpha D_2(t)|_{t=t_k} - ((u_2)_k - \gamma_2^\alpha (D_2)_k),
\end{cases} \quad (4.3)$$

$$\begin{cases} J_{11} &= r^\alpha - 2r^\alpha p^\alpha T_k - \xi_1^\alpha I_k + c_1^\alpha F_k - q_1^\alpha(D_1)_k, \quad J_{12} = -\xi_1^\alpha T_k, \\ J_{13} &= c_1^\alpha T_k, \quad J_{14} = -q_1^\alpha T_k, \quad J_{21} = \frac{2\rho^\alpha T_k I_k(h^\alpha + F_k^2)}{(h^\alpha + T_k^2 + F_k^2)^2} - \xi_2 I_k, \\ J_{22} &= \frac{\rho^\alpha T_k^2}{h^\alpha + T_k^2 + F_k^2} + \frac{\beta^\alpha (D_2)_k}{g^\alpha + (D_2)_k} - \xi_2^\alpha T_k - \mu^\alpha - q_2^\alpha(D_1)_k, \\ J_{23} &= -\frac{2\rho^\alpha T_k^2 I_k F_k}{(h^\alpha + T_k^2 + F_k^2)^2}, \quad J_{24} = -q_2^\alpha I_k, \quad J_{25} = \frac{\beta^\alpha g^\alpha I_k}{(g^\alpha + (D_2)_k)^2}, \\ J_{31} &= -c_2^\alpha F_k, \quad J_{33} = d^\alpha - 2d^\alpha \epsilon^\alpha F_k - c_2^\alpha T_k - q_3^\alpha(D_1)_k, \\ J_{34} &= -q_3^\alpha F_k, \quad J_{44} = -\gamma_1^\alpha, \quad J_{55} = -\gamma_2^\alpha. \end{cases} \quad (4.4)$$

## 4.2 Discrete adjoint equation

The discrete adjoint equation can be derived similarly as

$$\begin{aligned} {}_t^C\mathbb{D}_{t_f}^\alpha \varphi(t)|_{t=t_k} &= -\frac{1}{\Gamma(1-\alpha)} \sum_{j=k+1}^{N-1} \frac{\varphi_j - \varphi_{j-1}}{\Delta t} \int_{t_{j-1}}^{t_j} (s-t_k)^{-\alpha} ds \\ &= C_0 \sum_{j=k+1}^{N-1} (\varphi_{j-1} - \varphi_j)\left(\zeta_{j,k}^C\right), \end{aligned} \quad (4.5)$$

where $C_0 = -B_0$, $\zeta_{j,k}^C = ((j-k)^{1-\alpha} - (j-k-1)^{1-\alpha})$.

The discrete scheme (4.5) is applied to (2.8) to obtain the following system of equations for $N - 2 \geq k \geq 0$:

$$\begin{pmatrix} \delta_{k,k+1}^C C_0 - J_{11} & -J_{21} & -J_{31} & 0 & 0 \\ -J_{12} & \delta_{k,k+1}^C C_0 - J_{22} & 0 & 0 & 0 \\ -J_{13} & -J_{23} & \delta_{k,k+1}^C C_0 - J_{33} & 0 & 0 \\ -J_{14} & -J_{24} & -J_{34} & \delta_{k,k+1}^C C_0 - J_{44} & 0 \\ 0 & -J_{25} & 0 & 0 & \delta_{k,k+1}^C C_0 - J_{55} \end{pmatrix} \begin{pmatrix} (\lambda_1)_{k+1} \\ (\lambda_2)_{k+1} \\ (\lambda_3)_{k+1} \\ (\lambda_4)_{k+1} \\ (\lambda_5)_{k+1} \end{pmatrix}$$

$$= \begin{pmatrix} 1 + E_0 \sum_{j=k+2}^N (\lambda_1)_j \left(\delta_{k,j-1}^C - \delta_{k,j}^C\right) \\ E_0 \sum_{j=k+2}^N (\lambda_2)_j \left(\delta_{k,j-1}^C - \delta_{k,j}^C\right) \\ E_0 \sum_{j=k+2}^N (\lambda_3)_j \left(\delta_{k,j-1}^C - \delta_{k,j}^C\right) \\ E_0 \sum_{j=k+2}^N (\lambda_4)_j \left(\delta_{k,j-1}^C - \delta_{k,j}^C\right) \\ E_0 \sum_{j=k+2}^N (\lambda_5)_j \left(\delta_{k,j-1}^C - \delta_{k,j}^C\right) \end{pmatrix}. \quad (4.6)$$

## 5 Numerical results

In this section, we present some numerical results to compare three different treatment strategies: chemotherapy, immunotherapy and their combination. In addition, we investigate the effect of the order of differentiation $\alpha$ and the value of the decay rate of amount of chemotherapeutic drug $\gamma_1$ to the value of the cost functional $J_{\gamma_1}$. We set the values of the parameters as in Table 1, otherwise stated.

Table 1: Values of the parameters

| Parameter | Description | Value |
|---|---|---|
| $r$ | Per capita growth rate of tumor cells | 0.00431 |
| $p$ | Reciprocal carrying capacity of tumor cells | $1.02 \times 10^{-9}$ |
| $\xi_1$ | Competition term of tumor cells with immune cells | $6.41 \times 10^{-11}$ |
| $\xi_2$ | Competition term of immune cells with tumor cells | $3.42 \times 10^{-6}$ |
| $c_1$ | Competition term of fat cells with tumor cells | $\alpha_1$ |
| $c_2$ | Competition term of fat cells with tumor cells | $\alpha_2$ |
| $q_1$ | Response of tumor cells to chemotherapeutic drug | 0.08 |
| $q_2$ | Response of immune cells to immunotherapeutic drug | $2 \times 10^{-11}$ |
| $q_3$ | Response of fat cells to immunotherapeutic drug | $q_2$ |
| $s$ | Immune source rate | 0.33 |
| $\rho$ | Recruitment rate of immune cells by tumor cells | 0.0125 |
| $h$ | Immune response stimulated by tumor cells | 20,200,000 |
| $\mu$ | Per capita death rate of immune cells | 0.204 |
| $\beta$ | Recruitment rate of immune cells by $D_2$ | 0.125 |
| $g$ | Steepness coefficient of the $\beta$ | 20,000,000 |
| $d$ | Per capita growth rate of fat cells | $r/100$ |
| $\epsilon$ | Reciprocal carrying capacity of fat cells | $p$ |
| $\gamma_1$ | Decay rate of $D_1$ | 0.1 |
| $\gamma_2$ | Decay rate of $D_2$ | 1 |
| $\omega_1$ | Weight constant | 1 |
| $\omega_2$ | Weight constant | 2 |
| $T_0$ | Initial number of tumor cells | 2 |
| $I_0$ | Initial number of immune cells | .1 |
| $F_0$ | Initial number of fat cells | 1 |
| $D_{10}$ | Initial amount of chemotherapeutic drug | .5 |
| $D_{20}$ | Initial amount of immunotherapeutic drug | .5 |
| $u_1$ | Initial dose of $D_1$ | 0.5 |
| $u_2$ | Initial dose of $D_2$ | 0.5 |
| $t_f$ | Final time | 120 days |
| $\Delta t$ | Length of the time step | 0.25 |

We implement the so-called forward-backward sweep method under MATLAB to solve the OCP (2.6-2.7) [17, Chap. 5]. The algorithm can be summarized as follows:

Before investigating the optimal treatment strategy, we solve the uncontrolled cancer-obesity model for different values of the order of differentiation $\alpha$. In Fig. 1, we present the number of tumor cells, immune cells and fat cells obtained by taking $u_1 = 0$ and $u_2 = 0$ in the FDE (2.7). We observe that the tumor and fat cell population are increasing over time for different values of $\alpha$, while the number of immune cells approaches to a fixed value.

**Algorithm 1** Forward-backward sweep method
─────────────────────────────────────────────
1: Fix $\psi = -1$, $\delta = 0.001$.
2: Initiate the control $(u_{old})_1$, $(u_{old})_2$, the state $x_{old} = \{T_{old}, I_{old}, F_{old}, (D_{old})_1, (D_{old})_2\}$ and adjoint $p_{old} = \{(\lambda_{old})_1, (\lambda_{old})_2, (\lambda_{old})_3, (\lambda_{old})_4, (\lambda_{old})_5\}$.
3: **while** $\psi < 0$ **do**
4:    Solve the state equation (2.7) for $x = \{T, I, F, D_1, D_2\}$ using $T_0$, $I_0$, $F_0$, $D_{10}$, $D_{20}$, $(u_{old})_1$, $(u_{old})_2$ forward in time.
5:    Solve the adjoint equation (2.8) for $p = \{\lambda_1, \lambda_2, \lambda_3, \lambda_4, \lambda_5\}$ using $\lambda_i(t_f) = 0$, $x = \{T, I, F, D_1, D_2\}$ backward in time.
6:    Update the control using the gradient equation (2.10) to reach $u_1$, $u_2$.
7:    Compute $\chi_i = \delta\|x_i\| - \|x_i - (x_{old})_i\|$, $v_j = \delta\|u_j\| - \|u_j - (u_{old})_j\|$, $\rho_i = \delta\|p_i\| - \|p_i - (p_{old})_i\|$ and calculate $\psi = \min\{\chi_i, v_j, \rho_i\}$ for $i, j \in \{1, 2, 3, 4, 5\}, k \in \{1, 2\}$.
8: **end while**
─────────────────────────────────────────────

According to these results, an optimal treatment strategy is required to eliminate tumor burden as time passes.

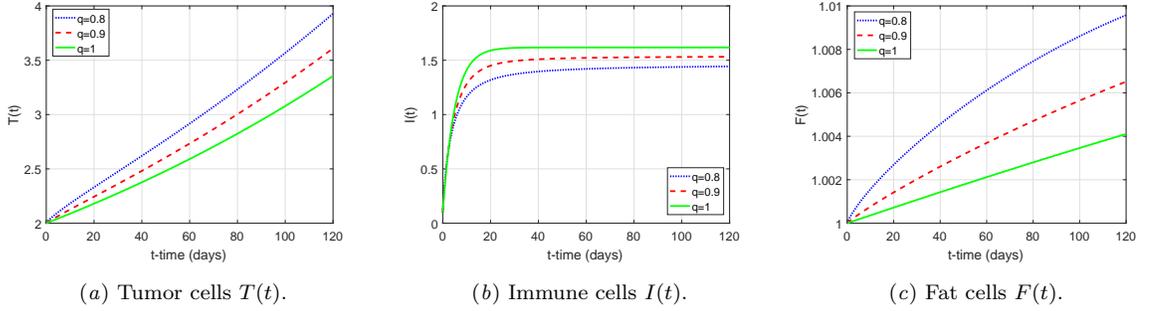

(a) Tumor cells $T(t)$.  (b) Immune cells $I(t)$.  (c) Fat cells $F(t)$.

Figure 1: Number of cells for uncontrolled case

## 5.1 Immunotherapeutic treatment

Firstly, we will examine the contribution of immunotherapy to cure the disease by taking $u_1 = 0$. In Fig. 2, the number of tumor cells, immune cells and fat cells are depicted. We can deduce from the figures that immunotherapeutic treatment has no positive effect on tumor burden, since the tumor cell population cannot be decreased.

In Fig. 3, we depict the immunotherapeutic drug concentration $D_2(t)$ and drug dose $u_2(t)$. We notice that immunotherapy has no contribution on tumor population although drug concentration is increased over time. In addition, the drug dose is quite small which is not enough to minimize the tumor population.

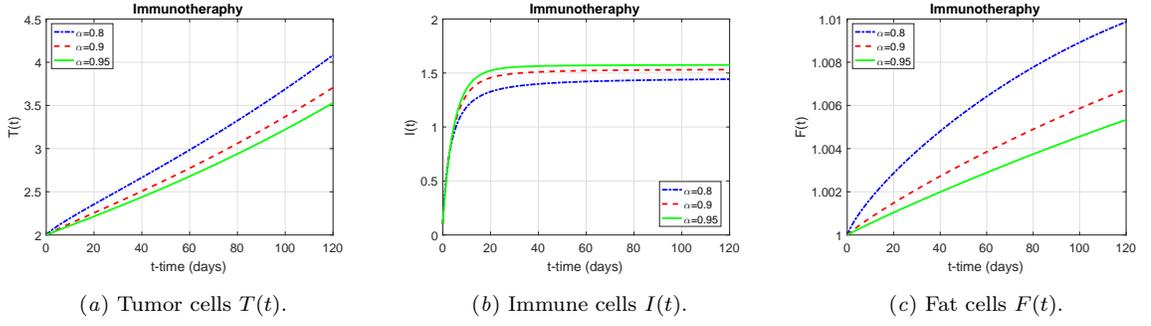

(a) Tumor cells $T(t)$.   (b) Immune cells $I(t)$.   (c) Fat cells $F(t)$.

Figure 2: Number of cells for immunotherapy

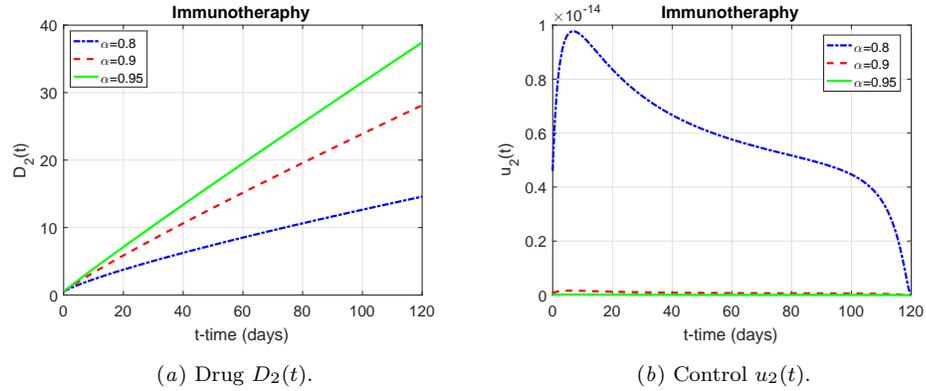

(a) Drug $D_2(t)$.   (b) Control $u_2(t)$.

Figure 3: Amount & dose of immunotherapeutic drug

## 5.2  Chemotherapeutic treatment

We proceed with chemotherapeutic treatment where $u_2 = 0$. In Fig. 4, tumor, immune and fat population are presented. We immediately see that tumor population is successfully minimized over time despite of increasing number of fat cells in the system, while the immune population is not destroyed.

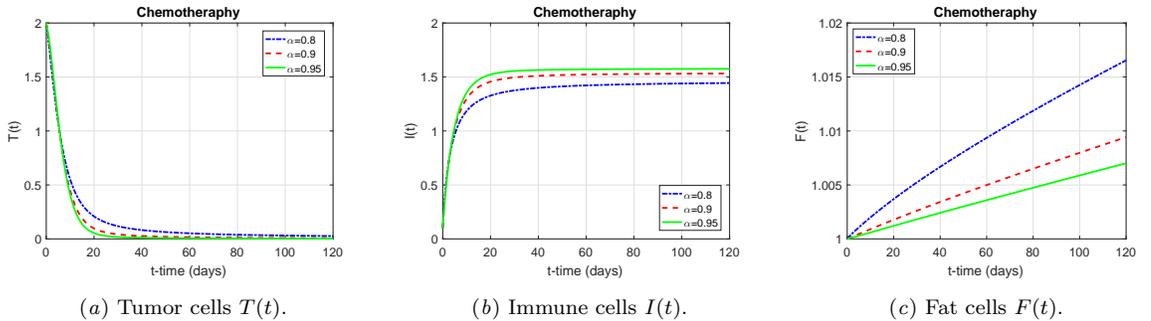

(a) Tumor cells $T(t)$.   (b) Immune cells $I(t)$.   (c) Fat cells $F(t)$.

Figure 4: Number of cells for chemotherapy

We investigate the effect of the order $\alpha$ and the decay rate $\gamma_1$ to the value of the cost functional $J_{\gamma_1}$ and we present the results in Table 2. Firstly, we notice that the value of $J_{\gamma_1}$ decreases as we increase $\alpha$. Moreover, a smaller value of $J_{\gamma_1}$ is measured for smaller values of the decay rate $\gamma_1$, since more drug is contained in the system. A FOCP with a smaller value of $\alpha$ can be regarded as a model containing information about history, so the system with a small value of $\alpha$ is more vulnerable to environmental changes.

Table 2: Values of the cost functional $J_{\gamma_1}$ for chemotherapy

| $\alpha$ | $J_{\gamma_1=0.1}$ | $J_{\gamma_1=0.5}$ | $J_{\gamma_1=0.9}$ |
|---|---|---|---|
| 0.8 | 88.0784 | 135.3614 | 177.9395 |
| 0.9 | 81.3347 | 117.0319 | 156.5642 |
| 0.95 | 79.3366 | 110.2287 | 147.9902 |

In Fig. 5, we present the chemotherapeutic drug concentration $D_1(t)$ and drug dose $u_1(t)$. We observe that the reverse relation between the order $\alpha$ and the value of $J_{\gamma_1}$ is revealed, since drug concentration gets higher as $\alpha$ increases. Thus, more tumor cells are destroyed due to higher drug concentration. On the other hand, the drug dose decreases as time passes, which means that the functional $J_{\gamma_1}$ has been minimized successfully over time.

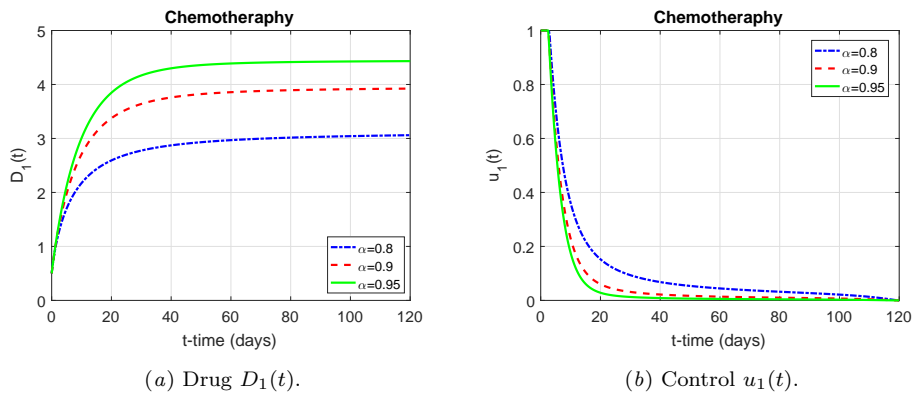

(a) Drug $D_1(t)$.

(b) Control $u_1(t)$.

Figure 5: Amount & dose of chemotherapeutic drug

### 5.3 Combination of chemotherapeutic & immunotherapeutic treatment

We proceed with the combination of chemotherapeutic and immunotherapeutic treatment where $u_1 \neq 0$ and $u_2 \neq 0$. In Fig. 6, the number of tumor cells, immune cells and fat cells are depicted. We observe that the tumor population is decreased. There is not a visible difference between chemotherapeutic treatment and combined therapy. Therefore, we will compute $J_{\gamma_1}$ to assess the contribution of combined therapy.

In Table 3, we present the values of $J_{\gamma_1}$ for different values of $\alpha$ and $\gamma_1$. As we increase $\gamma_1$, $J_{\gamma_1}$ decreases; while the fractional order has a reverse effect. When we compare Table 3 with Table 2, we deduce that combined therapy gives much better results than chemotherapy.

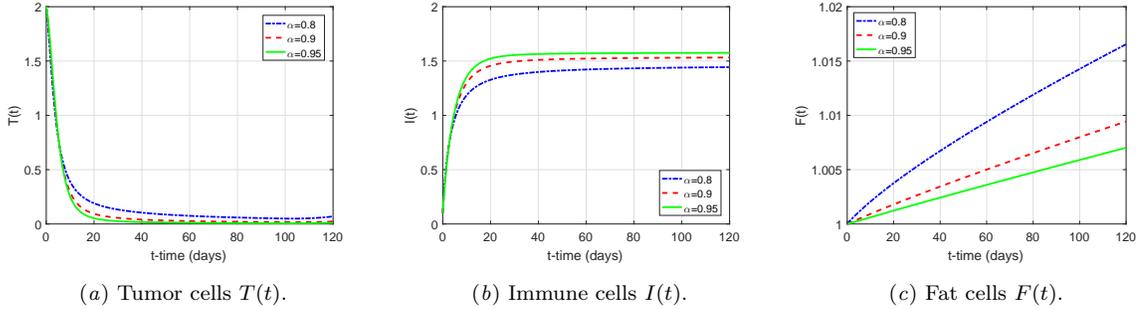

(a) Tumor cells $T(t)$.     (b) Immune cells $I(t)$.     (c) Fat cells $F(t)$.

Figure 6: Number of cells for combined chemotherapy & immunotherapy

On account of a stronger immune system, immune cells fight with tumor cells much better and tumor population shrinks in size. We think that a system including memory effect, which corresponds to the use of fractional order derivatives in the underlying model, can be regarded as a model associated with a population which can correctly adjust itself to environmental changes and which can rebound itself on account of previous experiences.

Table 3: Values of the cost functional $J_{\gamma_1}$ for chemotherapy & immunotherapy

| $\alpha$ | $J_{\gamma_1=0.1}$ | $J_{\gamma_1=0.5}$ | $J_{\gamma_1=0.9}$ |
|---|---|---|---|
| 0.8 | 36.7257 | 70.6701 | 95.7243 |
| 0.9 | 28.1606 | 53.7762 | 77.6948 |
| 0.95 | 24.9350 | 45.7107 | 69.0309 |

In Fig. 7, we present the drug concentration and dose for chemotherapeutic and immunotherapeutic treatment. As time passes, drug concentrations $D_1(t)$ and $D_2(t)$ decreases. Otherwise, it might be harmful for healthy and immune cells. Dose of drugs $u_1(t)$ and $u_2(t)$ lies within the admissible set and they decrease over time. Although the control $u_2$ seems too small, immunotherapeutic treatment leads $J_{\gamma_1}$ to decrease more than 50%.

## 6 Summary and conclusion

In this study, we investigate the optimal treatment strategy for a cancer-obesity model. Interactions between tumor cells, immune cells, fat cells, chemotherapeutic and immunotherapeutic drug concentrations are modeled with Caputo time fractional derivative. The aim is to find a pair of controls, which correspond to drug dose, to minimize the number of tumor cells over a finite time period, while finding the minimum drug dose injected to the system. We compare immunotherapy, chemotherapy and their combination to treat the tumor. We find out that mixed immune-chemotherapy gives the smallest value of the cost-functional $J$. Moreover, we notice that the values of $J$ decreases as we increase the order of differentia-tion $\alpha$. We think that a system including memory effect, which corresponds to the use of fractional order derivatives in the underlying model, can be regarded as a model associated with a population which can correctly adjust itself to environmental changes and which can rebound itself on account of previous experiences. As a future work, we will investigate the optimal treatment strategy through the use a fractional derivative with nonsingular kernel.

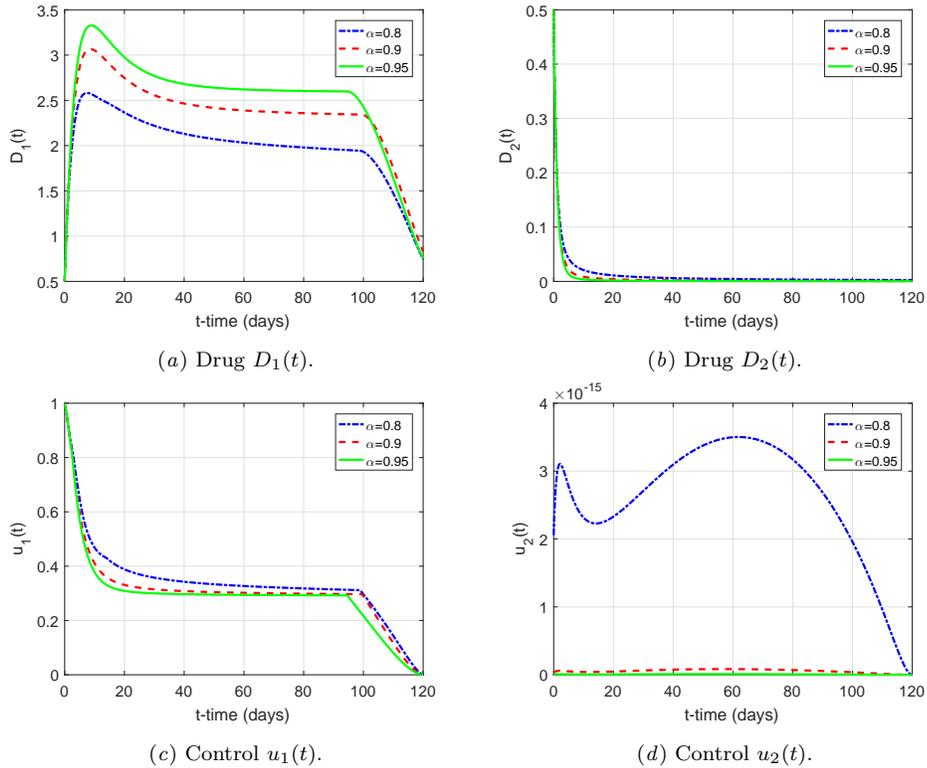

Figure 7: Amount of drug and drug control